\documentclass[preprint,12pt]{elsarticle}

\usepackage{amssymb}
\usepackage{amsmath}
\usepackage{bm}
\usepackage{float} 
\usepackage{algorithm}
\usepackage{algpseudocode}
\usepackage{color}
\usepackage{makecell}
\usepackage{multirow}
\usepackage{booktabs}
\usepackage{diagbox}
\usepackage{caption}
\usepackage{tabularx}
\usepackage{hyperref}
\usepackage{subcaption}

\newtheorem{remark}{Remark}[section]

\journal{arxiv}

\begin{document}
\begin{highlights}
	\item A domain decomposition method for the Finite Element simulation of surface acoustic wave (SAW) devices, with iteration numbers stable with respect to number of subdomains, frequency and mesh size.
	\item The proposed method exhibits excellent parallel performance and high computational efficiency. It maintains rapid convergence even for problems with millions of degrees of freedom, and yields a relative residual as low as \(10^{-13}\), which is comparable to that of direct solvers.
	\item A novel coarse space construction strategy is developed based on the intrinsic physical characteristics of SAWs, which captures the physical nature of SAW governing equations. This strategy is versatile and can be readily extended to the solver design for various types of SAW devices.
\end{highlights}

\begin{frontmatter}



	\title{A Novel Two-Level Domain Decomposition Method with Physically Motivated Coarse Space for Surface Acoustic Wave Device Simulations}


	\author[addr1]{Likun Qiu}
	\ead{qiulkmath@tongji.edu.cn}
	\author[addr2]{Jianguo Huang\corref{cor2}}
	\ead{jghuang@sjtu.edu.cn}
	\author[addr1]{Xuejun Xu\corref{cor1}}
	\ead{xuxj@tongji.edu.cn}
	\cortext[cor1]{Corresponding author: Xuejun Xu, E-mail: xuxj@tongji.edu.cn}
	\cortext[]{The work of the second author was partially supported by NSFC (Grant No. 12571390), and that of the third author was partially supported by NSFC (Grant No. 12331015).}

	\affiliation[addr1]{
		organization={School of Mathematical Sciences, Tongji University, and Key Laboratory of Intelligent Computing and Applications (Tongji University), Ministry of Education},
		city={Shanghai},
		postcode={200092},
		country={China}
	}
	\affiliation[addr2]{
		organization={School of Mathematical Sciences, and MOE-LSC, Shanghai Jiao Tong University},
		city={Shanghai},
		postcode={200240},
		country={China}
	}

	\begin{abstract}
		The surface acoustic wave (SAW) equations are fundamental in describing the operating principle of SAW devices, which are characterized by high operating frequency, multiphysics coupling, and defined on an unbounded domain. Finite element (FE) simulation is one of the most accurate numerical methods for such devices. However, the method usually leads to ill-conditioned linear systems with a large number of degrees of freedom (DOFs), limiting its application. To address this issue, based on the finite element tearing and interconnecting (FETI) algorithm, this paper proposes a domain decomposition method with iteration numbers stable with respect to number of subdomains, frequency and mesh size, named FETI-SAW. The main novelty is that the Dirichlet preconditioner combined with a Rayleigh-wave-based coarse space is constructed to produce a preconditioned iterative algorithm for the Lagrange multiplier linear system. In some sense, this method reflects the physical essence of the SAW equations. In numerical experiments, we take interdigital transducers (IDT) and SAW resonators as examples to demonstrate the efficiency of the FETI-SAW algorithm. Additionally, the MATLAB and C++ source codes accompanying this paper are provided.
	\end{abstract}


	%


\end{frontmatter}



\section{Introduction}\label{sect:Intro}
SAW devices are widely used in the areas of electronic and information engineering, such as radar, mobile communications, and the Internet of Things \cite{Saskatoon-2006,Zhao-2019,Mandal-Banerjee-2022}. With the rapid development of these fields, higher requirements have been put forward for the research and production of SAW devices. Especially in mobile communication systems, a large number of high-frequency, low-loss, high-performance SAW devices are needed \cite{Yang-2018}. This also imposes higher requirements on the efficiency and accuracy of the simulation methods. Phenomenological models such as the widely used equivalent-circuit model and coupled-mode theory can hardly accurately predict the loss in the passband of SAW filters. These methods often neglect the end effects of transducers and reflecting gratings, as well as the loss caused by bulk wave scattering \cite{Hashimoto-2025}. FE simulation is one of the most accurate simulation methods for SAW devices \cite{Saskatoon-2006,Ippolito-Kalantar-2003,Hofer-Finger-2006}. However,  the existing FE simulation has been hampered by large memory consumption and low computational speed. Therefore, it is necessary to develop high-performance solvers for the FE discretization of the SAW equations to achieve the efficient and high-precision FE simulations.

\begin{figure}[H]
	\centering
	\begin{minipage}[t]{0.45\textwidth}
		\centering
		\includegraphics[width=\textwidth]{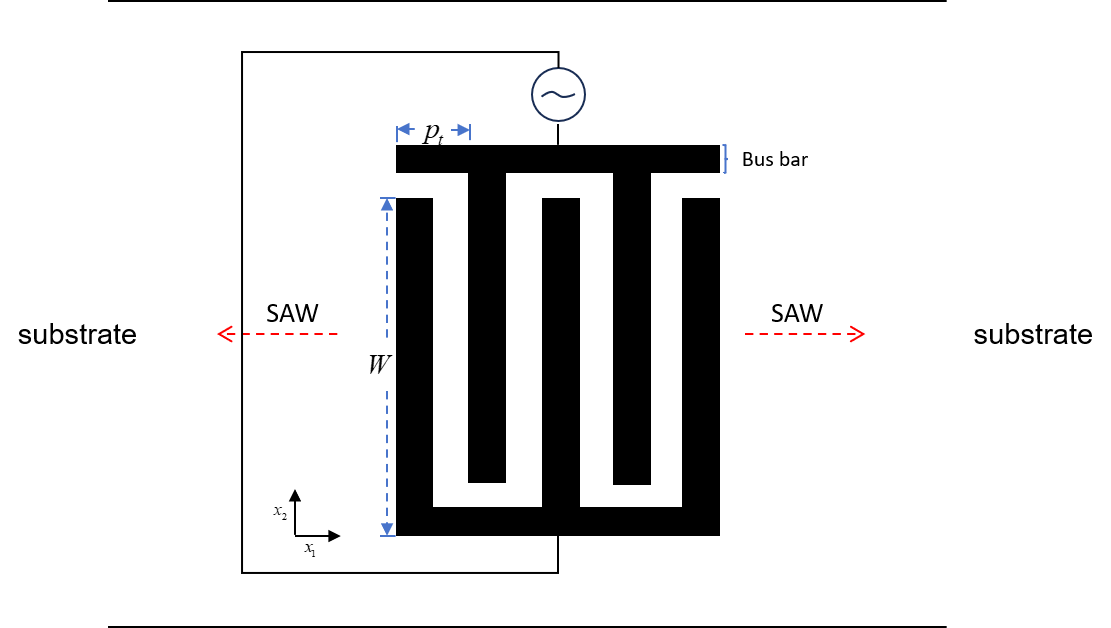}
	\end{minipage}
	\hfill
	\begin{minipage}[t]{0.45\textwidth}
		\centering
		\includegraphics[width=\textwidth]{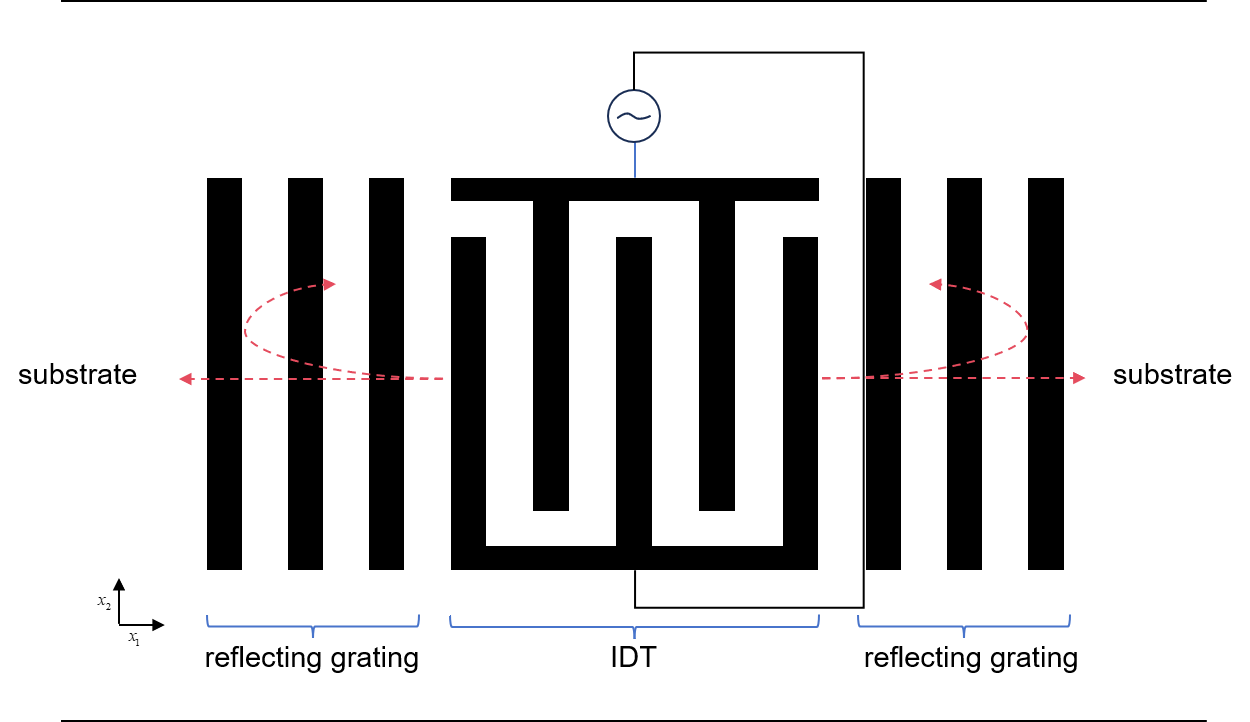}
	\end{minipage}
	\caption{Interdigital transducer (left); SAW resonator (right).}
	\label{fig:CommonSAWdevices}
\end{figure}

In SAW devices, IDT is the key structure. As shown in Figure \ref{fig:CommonSAWdevices} (left), the black regions are metallic thin-film electrodes and the blank region is a piezoelectric substrate. When an alternating current is applied to the metallic electrodes, the SAWs are excited on the surface of the piezoelectric substrates as a result of the piezoelectric effect of the substrates. The length \(p_t\) of two adjacent metal electrodes of the IDT is equal to half the typical wavelength of the SAWs.
In addition to transmitters, such structures as IDT can also serve as receivers and reflecting gratings in SAW devices \cite{Koskela-2001}. Figure \ref{fig:CommonSAWdevices} (right) shows a SAW resonator with an open-circuit reflecting grating. The SAW excited by the IDT superposes with the wave reflected back by the gratings, forming a standing wave on the surface of the piezoelectric substrate.
In SAW devices, the length \(W\) of the metal electrode is typically \(40\) to \(100\) times the typical wavelength of SAW. Therefore, the bus bar is usually assumed to only make the connected metal strips equipotential, the aperture direction is regarded as infinitely long, and all physical quantities are unchanged along the aperture direction (i.e., the \(x_2\) direction in Figure. \ref{fig:CommonSAWdevices}). Thus, the simulation and analysis of SAW devices are generally based on 2-D (c.f. \cite{Koskela-2018,Saskatoon-2006}) and 2.5-D models (c.f. \cite{El-Moussa-2009,Zhao-2019,Saskatoon-2006}). Compared with the 2-D model, the 2.5-D model reduces the problem size while accounting for the tangential strain and achieves higher accuracy \cite{Zhao-2019}. Therefore, in this paper, we mainly consider the 2.5-D model, and the result can be directly applied to the 2-D model.

In practice, we only need to simulate the region near the metal electrodes, whereas the piezoelectric substrate is much larger than this region, and hence the piezoelectric substrate is often assumed to be semi-infinite. Therefore, in the FE simulation of SAW devices, we first need to truncate the unbounded domain. To the best of our knowledge, there are two main approaches to addressing this wave problem: the PML method and the use of absorbing boundary conditions (ABC). Both of these approaches  truncate the unbounded domain by artificial boundary conditions. We refer the reader to \cite{ABC1,ABC2,ABC3,ABC4} for the ABC methods and \cite{PML1,Pled-Desceliers-2022,Leng-Ju-2019,Leng-Ju-2022} for the PML methods. Due to their simplicity and easy implementation, PML methods have been widely used in microacoustic devices including SAW and BAW devices (cf. \cite{Koskela-2018,PMLPiez1,PMLPiez3}). We will adopt this method in the paper as well.

As shown in Figure \ref{fig:CommonSAWdevices}, the substructures of SAW devices, such as the IDT and reflecting gratings, exhibit periodic structures. This characteristic of SAW devices provides an excellent entry point for designing high-performance solvers for the FE models. Based on the periodic structures of SAW devices, Koskela et al. proposed the hierarchical cascading algorithm (cf. \cite{Koskela-2018,Li-Koskela-Willemsen-Massey-2026}), which is an efficient algorithm for computing the response parameters of SAW devices by solving their FE models. However, the hierarchical tree of cascading operations in the hierarchical cascading algorithm severely limits its parallel performance, making it difficult to fully exploit the parallel computing capability of modern computers and thus inefficient for large-size problems. Domain decomposition methods are also well-suited for solving equations defined on domains with periodic structures, while respecting the memory hierarchy of modern parallel and distributed computing systems \cite{Toselli-Widlund-2005}.
In \cite{Qiu-2026}, the authors proposed a non-iterative domain decomposition method, referred to as the DirFETI algorithm, for the FE simulation of the IDT model. This approach is capable of solving problems involving more than million DOFs on a standard laptop, which is efficient in many aspects. Nevertheless, the DirFETI algorithm performs well only for IDT models with ideal periodic structures and cannot be extended to more complex SAW models such as SAW resonators. Furthermore, constrained by the use of the Sherman-Morrison-Woodbury formula for solving linear systems, the algorithm exhibits insufficient parallel capability, making it difficult to sufficiently exploit the advantages of high-performance computers. Therefore, in order to address the disadvantage of the previous method, this paper is intended to develop a novel domain decomposition solver for the FE simulation of more complex SAW devices to meet the engineering requirements.

The FETI algorithm, first proposed by Farhat and Roux \cite{Farhat-Roux-1991}, is a type of non-overlapping domain decomposition method. Combined with either the Dirichlet preconditioner or the lumped preconditioner, the FETI algorithm exhibits good scalability when solving symmetric positive definite elliptic partial differential equations \cite{Farhat-Mandel-Roux-1994,Maday-2007}. Wave problems are generally indefinite. To address this issue, Farhat et al. extended the FETI algorithm to FETI-H and FETI-DPH; these two algorithms demonstrate good scalability when solving the time-harmonic Helmholtz equations and the isotropic elastic wave equations, respectively (cf. \cite{Farhat-Macedo-Lesoinne-2000,Farhat-Li-Lesoinne-Avery-2005,Farhat-Macedo-Lesoinne-Roux-2000,Farhat-Chen-Risler-Roux-1998}). Subsequently, Jin et al. generalized the ideas of the FETI-H algorithm for solving the electromagnetic equations and proposed the FETI-DPEM algorithm in \cite{Li-Jin-2006,Li-Jin-2008}. The characteristics of SAWs share some similarities with the waves mentioned above. Motivated by the previous observation, we propose in this paper a high-performance solver called the FETI-SAW algorithm for the FE models of SAW devices, based on the FETI method.

The key to designing efficient FETI-type algorithms is to construct appropriate coarse spaces and preconditioners. Since wave propagation in SAW devices is complex and their FE models are severely ill-conditioned, it is ineffective to directly use the Dirichlet preconditioner or the plane wave coarse space introduced in \cite{Farhat-Li-Lesoinne-Avery-2005,Farhat-Macedo-Lesoinne-2000}. We notice that Rayleigh waves are the main component of the wave field in SAW devices \cite{Morgan-2007}. Motivated by this observation, we are tempted to construct a coarser space based on Rayleigh waves. On the other hand, SAW devices can be divided into an ordered arrangement of several types of subdomains. Using the Dirichlet preconditioner does not add excessive computational cost in each iteration, so we adopt this preconditioner in our algorithm. Based on the above ideas, we finally  propose a FETI-SAW algorithm, in which we only need to use one or two Rayleigh wave velocities to construct the coarse space. The algorithm therefore has a small sized coarse problem and excellent parallel performance, leading to high efficiency in implementation. In section \ref{sect:NumericalExample}, we employ some IDT and SAW resonator examples to demonstrate the efficiency and scalability of the FETI-SAW algorithm.

We end this section by introducing some notations frequently used later on. Throughout this paper, unless stated explicitly, we use Einstein's convention for summation whereby summation is implied when an index occurs exactly twice. For a vector with the $i$-th component $x_i$, we simply write $[x_i]$ or $\bm x$ as the vector itself. This convention also applies to any tensor. We denote the partial derivative of \(f\) with respect to \(x_i\) by \(\partial_i f\). We denote the transpose (not the conjugate transpose) of matrix \(\bm A\) as \(\bm A^{\rm T}\), and use \(\bm A^{*}\) to represent the conjugate transpose of \(\bm A\).

\section{Mathematical models}\label{sect:MathematicalModel}
\begin{figure}[H]
	\centering
	\includegraphics[width=\textwidth]{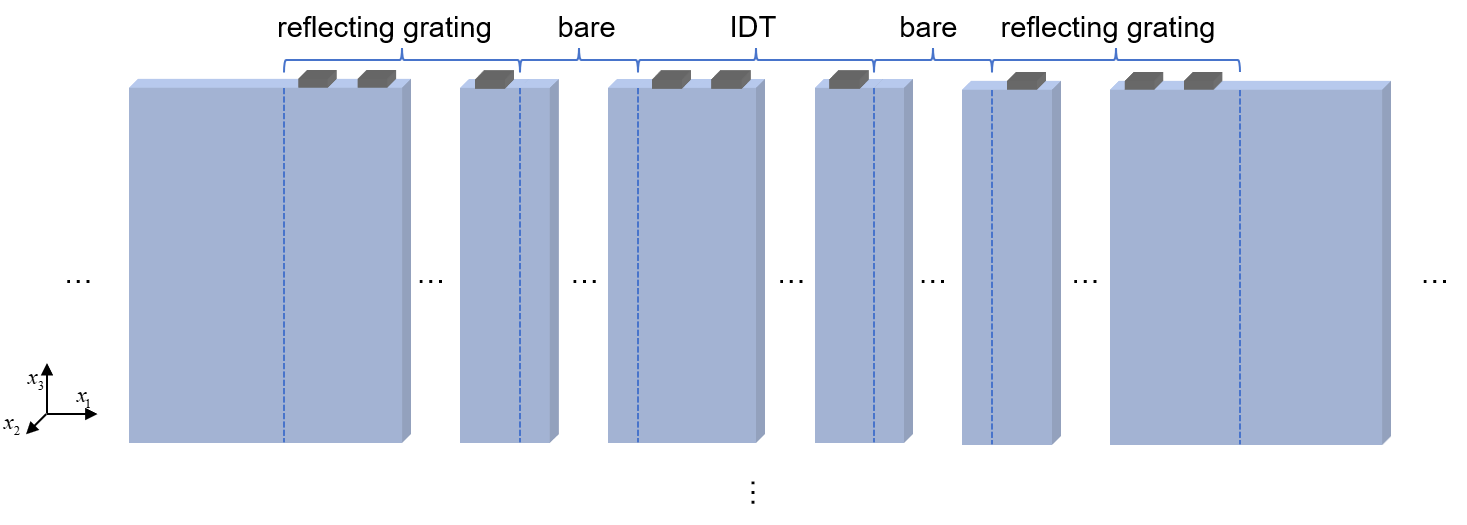}
	\caption{2.5-D model of SAW resonator. The grayish-blue region represent the semi-infinite piezoelectric substrate, and the black regions are the metal electrodes.}
	\label{fig:2_5Dresonator}
\end{figure}

\begin{figure}[H]
	\centering
	\includegraphics[width=0.75\textwidth]{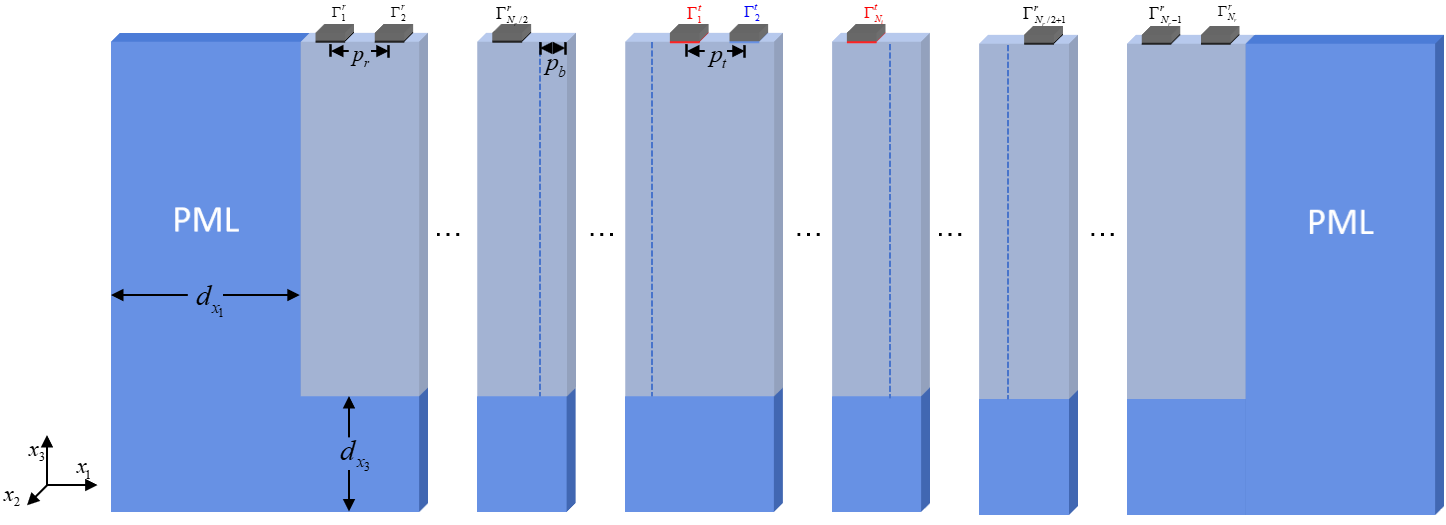}
	\caption{2.5-D model of SAW resonator truncated by PML.}
	\label{fig:2_5DresonatorPML}
\end{figure}

Without loss of generality, we take the SAW resonator in Figure \ref{fig:CommonSAWdevices} (right) as an example for illustration. As mentioned in Sect. \ref{sect:Intro}, we consider the 2.5-D model illustrated in Figure \ref{fig:2_5Dresonator}. In the 2.5-D model, the device dimension along the $x_2$-direction is significantly reduced.  When an alternating current is applied to the electrodes of the IDT, SAWs are excited on the surface of the piezoelectric substrate near the IDT, propagate along the $x_1$-direction, and decay exponentially along the negative $x_3$-direction. When these SAWs pass through the reflecting gratings, part of the wave is reflected, while the remainder continues to propagate along the surface.
In the simulation of SAW devices, the piezoelectric substrate is usually assumed to be semi-infinite, i.e., infinitely extended along the \(x_1\) and negative \(x_3\) directions. Therefore, we first adopt the PML method to truncate the unbounded domain, as shown in Figure \ref{fig:2_5DresonatorPML}.
For convenience, the \(x_1\)-coordinates of the left and right boundaries of the physical SAW resonator in Figure \ref{fig:2_5DresonatorPML} are denoted as \(x_l\) and \(x_r\), and the \(x_3\)-coordinate of the lower boundary as \(x_b\). Introduce
\begin{equation}\label{eqn:pmlDampingFun}\footnotesize
	\eta_1(x_1)=\left\{\begin{array}{ll}
		\eta_l(x_1) & \text{ in }[x_l-d_{x_1},x_l),\\
			              \eta_r(x_1)&\text{ in }(x_r,x_r+d_{x_1}], \\
		0           & \text{ otherwise,}
	\end{array}\right.\quad \eta_3(x_3)=\left\{\begin{array}{ll}
		\eta_b(x_3) & \text{ in }[x_b-d_{x_3},x_b), \\
		0           & \text{ otherwise,}
	\end{array}\right.
\end{equation}
and
\begin{equation}\label{eqn:dampingFunPML}
	\alpha_i=1-{\rm i}\eta_i(x_i)\quad(i=1,3),\quad \alpha_2=1,
\end{equation}
where \(d_{x_1}\) and \(d_{x_3}\) are the depths of the PML domain along \(x_1\) and \(x_3\),
\(\eta_l(x_1)\), \(\eta_r(x_1)\), \(\eta_b(x_3)\) are continuous positive-valued functions that increases monotonically with depth,
\(\mathrm{i}=\sqrt{-1}\) denotes the imaginary unit, and
\[
	\eta_l(x_1)=\eta_r(x_l+x_r-x_1).
\]
Readers interested in the principle and implementation details of the PML method may refer to \cite{Pled-Desceliers-2022}.

Let \(\Omega^p\) be the entire piezoelectric substrate, \(\Omega^e\) be the domain containing all the electrodes, \(\Omega^{\rm PML}\) be the PML domain and \(\Gamma_{\rm ext}=\{\bm x:x_1=x_l-d_{x_1}\text{ or }x_1=x_r+d_{x_1}\text{ or }x_3=x_b-d_{x_3}\}\). We denote \(\Gamma_i^r\) (\(i=1,\cdots,N_r\)) and \(\Gamma_i^t\) (\(i=1,\cdots,N_t\)) as the electrode-substrate interfaces in the reflecting gratings and IDTs, respectively. Metal electrodes are generally regarded as equipotential bodies; therefore, we neglect the electrical physical quantities within \(\Omega^e\). Then, the governing equations defined in the entire domain are
\begin{equation}\label{eqn:SAWEqnTruncated}\footnotesize
	\left\{\begin{array}{ll}
		\tilde{\partial}_j({c}^p_{ijkl}\tilde{\partial}_ku_l+{e}_{kij}^p\tilde{\partial}_k\phi)+\rho^p\omega^2 u_i=0, & \text{in }\Omega^p\text{ and }\Omega^{\rm PML},                                \\
		\tilde{\partial}_i({e}_{ikl}^p\tilde{\partial}_ku_l-\varepsilon^p_{ik}\tilde{\partial}_k\phi)=0,              & \text{in }\Omega^p\text{ and }\Omega^{\rm PML},                                \\
		\partial_j{c}^e_{ijkl}\partial_ku_l+\rho^e\omega^2 u_i=0,                                                     & \text{in }\Omega^e,                                                            \\
		\partial_2u_i=0,\quad \partial_2\phi=0,                                                                       & \text{ in }\Omega^p\cup\Omega^{\rm PML},                                       \\
		\partial_2u_i=0,                                                                                              & \text{ in }\Omega^e,                                                           \\
		u_i=0,\quad \phi=0,                                                                                           & \text{on }\Gamma_{\rm ext},                                                    \\
		\phi=\phi_0,                                                                                                  & \text{on } \Gamma_{2k+1}^t, k=0,....,\left\lfloor\frac{N_t-1}{2}\right\rfloor, \\
		\phi=0,                                                                                                       & \text{on }\Gamma_{2k}^t, k=1,....,\left\lfloor\frac{N_t}{2}\right\rfloor,      \\
		\phi(\bm x)=\phi(\bm x'),                                                                                     & \bm x,\bm x'\in \Gamma_i^r\quad(i=1,\cdots, N_r),                              \\
		u_i|_{\Omega^p}=u_i|_{\Omega^e},                                                                              & \text{ on }\partial\Omega^p\cap\partial\Omega^e                                \\
		\int_{\Gamma_i^r}D_jn_j{\rm dS}=0,                                                                            & i=1,\cdots, N_r,                                                               \\
		\sigma_{ij}n_j=0,\quad D_in_i=0,                                                                              & \text{on the rest boundary face,}
	\end{array}\right.
\end{equation}
where \(\tilde{\partial}_i=\frac{1}{\alpha_i}\partial_i\), $\bm n=[n_i]$ is the unit outward normal, $\phi$ is the electric potential, $[c_{ijkl}^{\cdot}]$ is the elasticity tensor, $[e_{kij}^{\cdot}]$ is the piezoelectric tensor,  $[\varepsilon_{ik}^{\cdot}]$ is the dielectric permittivity tensor, $\bm \sigma=[\sigma_{ij}]$ is the stress tensor, ${\bm u}=[u_j]$ is the displacement field, \(\phi\) is the electric potential, \(\phi_0\) is the voltage applied to the metal electrode of IDTs, $\bm D=[D_j]$ is the electric displacement vector, $\omega$ is the angular frequency, $\rho^{\cdot}$ is the density, and the superscripts $p$ and $e$ represent the material parameters correspond to the piezoelectric substrates and the electrodes, respectively.

In the piezoelectric substrate and PML regions, the relationship between \(\sigma_{ij}\), \(D_i\), \(u_i\) and \(\phi\) can be described by
\[
	\left\{\begin{array}{l}
		\sigma_{ij}=c_{ijkl}^p\tilde{\partial}_ku_l+e_{kij}^p\tilde{\partial}_k\phi, \\
		D_i=e_{ikl}^p\tilde{\partial}_ku_l-\varepsilon_{ik}^p\tilde{\partial}_k\phi.
	\end{array}\right.
\]
In the electrodes, the relationship between \(\sigma_{ij}\) and \(u_i\) is
\[
	\sigma_{ij}=c_{ijkl}^e\partial_ku_l.
\]
\section{Dimensionless formulation}
In practical, the material parameters in the SAW problems vary significantly in magnitude: the order of magnitude of the elastic tensor is typically around \(10^{10}\), while that of the dielectric tensor is approximately \(10^{-12}\). It degrades the numerical stability of numerical algorithms. In engineering applications, dimensionlessization is a conventional method to reduce the impact of the magnitude differences in material parameters on the stability and accuracy of numerical simulations. In this subsection, we introduce the dimensionlessization method for equations \eqref{eqn:SAWEqnTruncated} in detail.

First, we select the characteristic quantities of displacement \((U_0)\), piezoelectric tensor \((e_0)\), geometric dimension \((L_0)\), density \((\rho_0)\), and elastic tensor \((c_0)\) as follows:
\begin{table}[H]
	\footnotesize
	\centering
	\begin{tabular}{|l|l|}
		\hline
		\(U_0\)    & empirical estimation                                                           \\ \hline
		\(e_0\)    & \(\max\{|e_{kij}|\}\)                                                          \\ \hline
		\(L_0\)    & minimum of length/width/height of all subdomains in Figure \ref{fig:subDomain} \\ \hline
		\(\rho_0\) & \(\max\{\rho^p,\rho^e\}\)                                                      \\ \hline
		\(c_0\)    & \(\max\{|c_{ijkl}^p|,|c_{ijkl}^e|\}\)                                          \\ \hline
	\end{tabular}
	\caption{Definitions of characteristic quantities for some quantities.}
	\label{tab:dimensionless_quantities}
\end{table}
Through the relationships among the quantities, the characteristic quantities of the other quantities can be derived, as shown in Table \ref{tab:dimensionless_other_quantities}.
\begin{table}[H]
	\footnotesize
	\centering
	\begin{tabular}{|l|l|}
		\hline
		permittivity characteristic quantity \(\varepsilon_0\) & \(\frac{e_0^2}{c_0}\)               \\ \hline
		electric potential characteristic quantity \(\Phi_0\)  & \(\frac{e_0U_0}{\epsilon_0}\)       \\ \hline
		stress characteristic quantity \(\sigma_0\)            & \(\frac{c_0U_0}{L_0}\)              \\ \hline
		electric displacement characteristic quantity \(D_0\)  & \(\frac{\epsilon_0\Phi_0}{L_0}\)    \\ \hline
		angular frequency characteristic quantity \(\omega_0\) & \(\sqrt{\frac{c_0}{\rho_0 L_0^2}}\) \\ \hline
	\end{tabular}
	\caption{Definitions of characteristic quantities for other quantities.}
	\label{tab:dimensionless_other_quantities}
\end{table}
For any quantity, the relationship between its original form \(\mathcal{P}\), characteristic quantity \(\mathcal{P}_0\), and dimensionless form \(\bar{\mathcal{P}}\) is :
\begin{equation}\label{eqn:dimensionlessRelationship}
	\mathcal{P}=\mathcal{P}_0\bar{\mathcal{P}}.
\end{equation}

By substituting all physical quantities and material parameters in \eqref{eqn:SAWEqnTruncated} with their dimensionless forms by the relationship \eqref{eqn:dimensionlessRelationship}, we obtain the dimensionless equations. In the subsequent sections of this paper, we directly study the original equations for simplicity, and the techniques developed later on   can be extended to dimensionless problems in a straightforward way.
\section{FETI-SAW algorithm}
\subsection{Basic framework of the FETI algorithm}
\begin{figure}[H]
	\centering
	\includegraphics[width=0.8\textwidth]{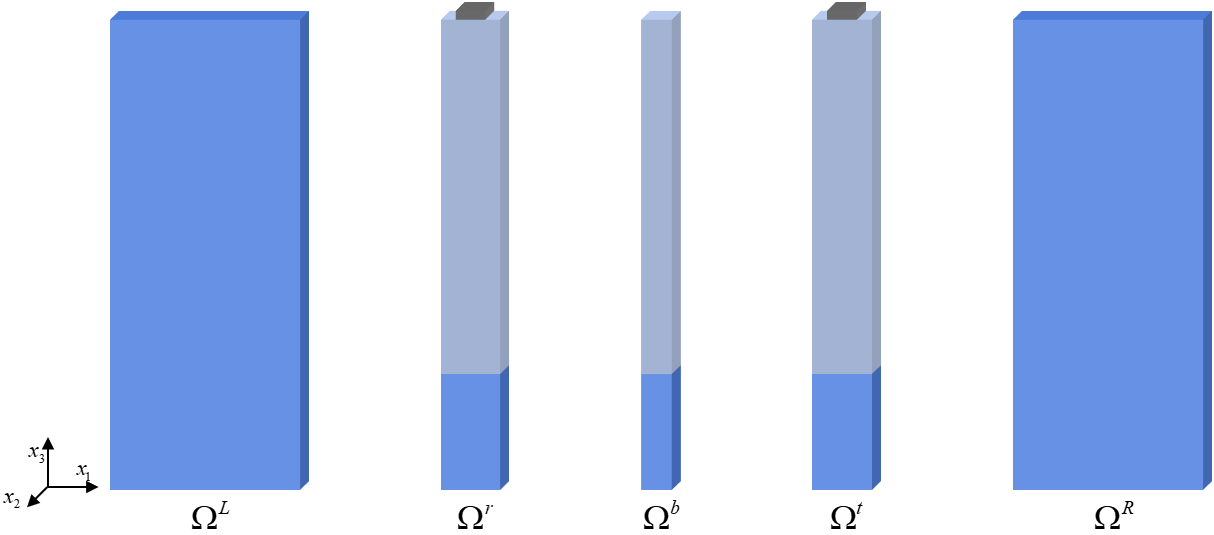}
	\caption{Five types of subdomain.}
	\label{fig:subDomain}
\end{figure}

As illustrated in Figure \ref{fig:2_5DresonatorPML}, the entire computational domain can be divided into five types of subdomains, as depicted in Figure \ref{fig:subDomain}. From left to right, these subdomains are as follows: the left PML domain \(\Omega^L\), the reflecting grating substructure \(\Omega^r\), the bare substrate  substructure \(\Omega^b\), the IDT substructure \(\Omega^t\), and the right PML domain \(\Omega^R\). 

We use ``width'' to denote the length of each subdomain along the \(x_1\) direction. Then, according to Figure \ref{fig:2_5DresonatorPML}, the widths of each subdomain from left to right are as follows: \(d_{x_1}\), \(p_r\), \(p_b\), \(p_t\) and \(d_{x_1}\). The domain illustrated in Figure \ref{fig:2_5DresonatorPML} can be represented as a sequence of these five types of subdomains:
\begin{equation}\label{eqn:subDomainArray}
	\Omega^L,\underbrace{\Omega^r,\cdots, \Omega^r}_{N_r/2},\underbrace{\Omega^b,\cdots, \Omega^b}_{N_b/2},\underbrace{\Omega^t,\cdots, \Omega^t}_{N_t},\underbrace{\Omega^b,\cdots, \Omega^b}_{N_b/2},\underbrace{\Omega^r,\cdots, \Omega^r}_{N_r/2},\Omega^R,
\end{equation}
where \(N_r\), \(N_b\), and \(N_t\) denote the number of substructures contained in the reflecting gratings, the bare substrate, and the IDT, respectively.

We discretize equations \eqref{eqn:SAWEqnTruncated} by the Lagrange FEM, ensuring uniform meshing for subdomains of the same type and matching grid points at the interfaces between adjacent subdomains. In this way, the matrices of the equations restricted to the same type subdomains are identical and only the right-hands of the odd-indexed IDT substructures  are nonzero vectors. For convenience, we use \(\Omega^{(m)}\) \((m=1,\cdots, N)\) to denote the \(m\)-th subdomain of the subdomain sequence \eqref{eqn:subDomainArray}, where \(N=2+N_b+N_t+N_r\) be the total number of the subdomains.

We denote the discrete form of \eqref{eqn:SAWEqnTruncated} as
\begin{equation}\label{eqn:originalEqn}
	\bm K\bm X=\bm f,
\end{equation}
where $\bm X$ consists of the numerical solutions for both displacement and electric potential. The matrices, right-hand side vectors and numerical solutions restricted to each subdomain are then denoted as $\bm K^{(m)}$, $\bm f^{(m)}$ and $\bm X^{(m)}$, respectively.
Like the classic FETI-type algorithms, we introduce signed boolean matrices \(\bm B^{(m)}\) that extracts from a given subdomain vector a signed interface boundary component. The sign of that component is specified by the orientation of the normal to the corresponding subdomain interface. Using this notation, the continuity of the subdomain solutions across the subdomain interfaces can be written as
\[\sum_{m=1}^N\bm B^{(m)}\bm X^{(m)}=0.\]
Then, we define the Lagrangian function as
\begin{equation}\label{eqn:LagrangianFun}\footnotesize
	\mathcal{L}(\bm Y^{(1)},\cdots,\bm Y^{(N)},\bm\lambda)=\frac{1}{2}\sum_{m=1}^N {\bm Y^{(m)}}^{\rm T}\bm K^{(m)}\bm Y^{(m)}-\sum_{m=1}^N{\bm f^{(m)}}^{\rm T}\bm Y^{(m)}+\bm\lambda^{\rm T}\sum_{m=1}^N\bm B^{(m)}\bm Y^{(m)},
\end{equation}
where \(\bm\lambda\) is the Lagrange multipliter defined on the interfaces.

Solving the original equation \eqref{eqn:originalEqn} is equivalent to find the stationary points of the Lagrangian function \eqref{eqn:LagrangianFun}. The Euler equations associated with the Lagrangian  function are given by
\begin{equation}\label{eqn:EulerEqn}
	\left\{\begin{array}{ll}
		\bm K^{(m)}\bm X^{(m)}=\bm f^{(m)}-{\bm B^{(m)}}^{\rm T}\bm\lambda, & m=1,\cdots,N \\
		\sum_{m=1}^N\bm B^{(m)}\bm X^{(m)}=\bm 0.
	\end{array}\right.
\end{equation}
Eliminating \(\bm X^{(m)}\) from equation \eqref{eqn:EulerEqn}, we obtain the following linear system for \(\bm\lambda\).
\begin{equation}\label{eqn:LagMultipliersEqn}
	\bm F\bm\lambda=\bm d,
\end{equation}
where
\[
	\bm F=\sum_{m=1}^N\bm B^{(m)}{\bm K^{(m)}}^{-1}{\bm B^{(m)}}^{\rm T},\quad \bm d=\sum_{m=1}^N\bm B^{(m)}{\bm K^{(m)}}^{-1}\bm f^{(m)}.
\]
After obtaining the Lagrange multiplier by solving equation \eqref{eqn:LagMultipliersEqn}, we can compute the numerical solutions of the original equation by
\begin{equation}\label{eqn:computeSourceSol}
	\bm X^{(m)}={\bm K^{(m)}}^{-1}\bm f^{(m)}-{\bm K^{(m)}}^{-1}{\bm B^{(m)}}^{\rm T}\bm\lambda.
\end{equation}
The key to design an efficient FETI-type algorithm is to develop an efficient solver for equation \eqref{eqn:LagMultipliersEqn}. In the next two subsections, we introduce the classic preconditioner of the FETI algorithm and how to design a two-level algorithm based on the Rayleigh wave.
\subsection{Preconditioners}
In this subsection, we introduce two classic preconditioners of the FETI algorithm. We use subscripts \(\Gamma\) and \(i\) to denote symbols associated with the interface and the remaining DOFs, respectively. Then, the Dirichlet preconditioner and the lumped preconditioner are given by:
\begin{enumerate}[(1)]
	\item Dirichlet preconditioner\\
	      \begin{equation}\label{eqn:DiriPrecond}
		      \bm M_D^{-1}=\sum_{m=1}^N\bm B^{(m)}\left[\begin{array}{cc}
				      \bm0 & \bm0                       \\
				      \bm0 & \bm S_{\Gamma\Gamma}^{(m)}
			      \end{array}\right]{\bm B^{(m)}}^{\rm T},
	      \end{equation}
	      where
	      \[
		      \bm S_{\Gamma\Gamma}^{(m)}=\bm K_{\Gamma\Gamma}^{(m)}-\bm K_{\Gamma i}^{(m)}{\bm K_{i i}^{(m)}}^{-1}\bm K_{i\Gamma}^{(m)};
	      \]
	\item lumped preconditioned\\
	      \begin{equation}\label{eqn:LumpedPrecond}
		      \bm M_L^{-1}=\sum_{m=1}^N\bm B^{(m)}\left[\begin{array}{cc}
				      \bm0 & \bm0                       \\
				      \bm0 & \bm K_{\Gamma\Gamma}^{(m)}
			      \end{array}\right]{\bm B^{(m)}}^{\rm T}.
	      \end{equation}
\end{enumerate}

The Dirichlet preconditioner is mathematically optimal \cite{Farhat-Mandel-Roux-1994}, but computing the Schur complement  \(\bm S_{\Gamma\Gamma}\) incurs a considerable computational cost at each iterative step. The lumped preconditioner is an approximation of the Dirichlet preconditioner, which significantly reduces the computational cost per iterative step, but its reduction in the iteration numbers is rather limited compared with the Dirichlet preconditioner. The roots and mechanical interpretation of these two preconditioners can be found in \cite{Farhat-Mandel-Roux-1994,Farhat-Roux-1991}.

However, when solving the SAW equations, using either the Dirichlet preconditioner or the lumped preconditioner alone yields a very limited reduction in the number of iterations. Even worse, the additional computational cost introduced in each iterative step can even lead to a substantial increase in the total runtime of the algorithm.

\subsection{Construction of the Rayleigh wave coarse space}
A key contribution of this paper is to construct the coarse space based on the propagation characteristics of SAW, which can produce an FETI-SAW algorithm with very good scalability.

Let \(\bm r_k\) and \(\bm\lambda_k\) denote the residual and solution of \eqref{eqn:LagMultipliersEqn} at the \(k\)-th iteration step, where
\begin{equation}\label{eqn:residual}
	\bm r_k=\bm d-\bm F\bm\lambda_k.
\end{equation}
From \eqref{eqn:residual} and the definition of \(\bm d\) and \(\bm F\), we can deduce that
\[
	\bm r_k=\sum_{m=1}^N\bm B^{(m)}\bm X^{(m)}_k,
\]
which shows that \(\bm r_k\) is the jump of the original solution \(\bm X^{(m)}\) on the subdomain interfaces. Therefore, if we can accurately find a set of basis \(\bm Q\) of the original solution restricted to the interfaces and enforce \(\bm r_k\) orthogonal to \(\bm Q\) at each iteration \(k\), i.e.
\begin{equation}\label{eqn:Orth2Q}
	\bm Q^*\bm r_k=0,
\end{equation}
then the iterative method converges in one iteration. However, constructing such a matrix \(\bm Q\) is impractical. First, it is nearly impossible to find a set of basis that accurately characterize the characteristics of the wave field in SAW devices. Furthermore, the dimension of such a set of basis is extremely large, and the computational cost of \eqref{eqn:Orth2Q} is similar to solving equation \eqref{eqn:LagMultipliersEqn} using a direct method.

It is hard to precisely characterize the characteristics of the wave field in SAW devices, but we can approximately characterize the dominant wave modes. In SAW devices, the Rayleigh wave is generally the dominant wave mode. We therefore attempt to construct the coarse space based on this propagation characteristics of the Rayleigh wave.

In \cite{Morgan-2007} and \cite{Saskatoon-2006}, a method is given for computing Rayleigh wave solutions in semi-infinite piezoelectric substrate. According to the properties of Rayleigh waves, which propagate along the $x_1$-direction and decay exponentially along the negative $x_3$-direction, we assume that the time-harmonic displacement and electric field take the form:
\begin{equation}\label{eqn:partialSAW}
	\begin{aligned}
		u_k'  & = u_{k,0}'\exp\bigl(-\mathrm{i}( \beta x_1-\gamma x_3-\omega t)\bigr),  \\
		\phi' & = \phi_{0}'\exp\bigl(-\mathrm{i}( \beta x_1-\gamma x_3-\omega t)\bigr),
	\end{aligned}
\end{equation}
where $\beta$ is the wave number of SAW, $t$ denotes time, and $\gamma$, $u'_{k,0}$, $\phi'_0$ are unknowns.

Substituting \eqref{eqn:partialSAW} into the following fundamental piezoelectric wave equation.
\begin{equation}\label{eqn:basicPiezEqn}
	\left\{\begin{array}{l}
		\partial_j{c}_{ijkl}^p\partial_ku_l+{e}_{kij}^p\partial_k\phi-\rho^p \ddot{u}_i=0, \\
		\partial_i{e}_{ikl}^p\partial_ku_l-\varepsilon_{ik}^p\partial_k\phi=0,
	\end{array}\right.
\end{equation}
where \(\ddot{u}_i\) denotes the second-order time derivative. Then, we obtain a quadratic polynomial eigenvalue problem given by
\begin{equation}\label{eqn:QPEP}
	(\mu^2\bm A_2+\mu\bm A_1+\bm A_0)\bm X^\prime=0,
\end{equation}
where \(\mu=\gamma/\beta\), \(\bm X^\prime=[u_{1,0}^\prime,u_{2,0}^\prime,u_{3,0}^\prime,\phi_{0}^\prime]^{\rm T}\),

\[\footnotesize
	\begin{array}{ll}
		\bm A_2=- & \left[\begin{array}{llll}
				            c^p_{1331} & c^p_{1332} & c^p_{1333} & e^p_{313}           \\
				            c^p_{2331} & c^p_{2332} & c^p_{2333} & e^p_{323}           \\
				            c^p_{3331} & c^p_{3332} & c^p_{3333} & e^p_{333}           \\
				            e^p_{313}  & e^p_{323}  & e^p_{333}  & -\varepsilon^p_{33}
			            \end{array}\right],                                                       \\
		\bm A_1=  & \left[\begin{array}{llll}
				            c^p_{1131}+c^p_{1311} & c^p_{1132}+c^p_{1312} & c^p_{1133}+c^p_{1313} & e^p_{113}+e^p_{311}                      \\
				            c^p_{2131}+c^p_{2311} & c^p_{2132}+c^p_{2312} & c^p_{2133}+c^p_{2313} & e_{123}^p+e_{321}^p                      \\
				            c^p_{3131}+c^p_{3311} & c^p_{3132}+c^p_{3312} & c^p_{3133}+c^p_{3313} & e_{133}^p+e_{331}^p                      \\
				            e_{113}^p+e_{311}^p   & e_{123}^p+e_{321}^p   & e_{133}^p+e_{331}^p   & -(\varepsilon_{13}^p+\varepsilon_{31}^p)
			            \end{array}\right], \\
		\bm A_0=- & \left[\begin{array}{llll}
				            c^p_{1111}-\rho^p v^2 & c^p_{1112}            & c^p_{1113}            & e^p_{111}           \\
				            c^p_{2111}            & c^p_{2112}-\rho^p v^2 & c^p_{2113}            & e^p_{121}           \\
				            c^p_{3111}            & c^p_{3112}            & c^p_{3113}-\rho^p v^2 & e^p_{131}           \\
				            e^p_{111}             & e^p_{121}             & e^p_{131}             & -\varepsilon_{11}^p
			            \end{array}\right],
	\end{array}
\]
and \(v=\omega/\beta\) is the velocity of SAW.

Given a wave velocity \(v\), we can get eight sets of eigenpairs \((\mu,\bm X^\prime)\) by solving \eqref{eqn:QPEP}. These eight eigenpairs can be divided into four groups, corresponding to two quasi-shear wave terms, one quasi-longitudinal wave term, and one electromagnetic wave term, respectively \cite{Auld-1973}. When \(v\) is lower than the quasi-shear horizontal (Q-SH) wave velocity, these eight eigenvalues form four pairs of complex conjugates. Since the Rayleigh waves propagate along the \(x_1\) direction and decay exponentially along the negative \(x_3\) direction, we select the four eigenpairs with negative imaginary parts. Substituting these four eigenpairs into \eqref{eqn:partialSAW}, we obtain a set of Rayleigh wave basis functions corresponding to wave velocity \(v\).

Based on the SAW device structure and the material of the piezoelectric substrate, we select \(N_v\) velocities \(\{v_s\}_{s=1}^{N_v}\) and construct the corresponding SAW bases via the aforementioned method. Then, we denote the selected eigenvalues as \(\{\mu_{v_{s}}^j\}_{j=1}^{4}\). Let \(\bm P^{(m)}_+\) denote the set of all points on the interface of \(\Omega^{(i)}\) for which the unit outward normal vector is \((1,0,0)^{\rm T}\). We use \(u_{k,\mu_{v_{s}}^j}^\prime\) and \(\phi_{\mu_{v_{s}}^j}^\prime\) to denote the SAW basis related to  \(\mu_{v_{s}}^j\). Then, we define the submatrix \(\bm Q^{(m)}\) of \(\bm Q\) as
\begin{equation}\label{eqn:subCoarseMat}
	\bm Q^{(m)}=\left[\begin{array}{llll}
			\bm Q^{(m)}_{v_1} & \bm Q^{(m)}_{v_2} & \cdots & \bm Q^{(m)}_{v_{N_v}}
		\end{array}\right].
\end{equation}
where
\[
	\begin{aligned}
		 & \bm Q^{(m)}_{v_s}=\left[\begin{array}{llll}
				                     \bm q^{(m)}_{\mu_{v_{s}}^1} & \bm q^{(m)}_{\mu_{v_{s}}^2} & \bm q^{(m)}_{\mu_{v_{s}}^3} & \bm q^{(m)}_{\mu_{v_{s}}^4}
			                     \end{array}\right],                                                    \\
		 & \bm q^{(m)}_{\mu_{v_{s}}^j}=\left[\begin{array}{llll}
				                               u_{1,\mu_{v_{s}}^j}^\prime(P^{(m)}_+) & u_{2,\mu_{v_{s}}^j}^\prime(P^{(m)}_+) & u_{3,\mu_{v_{s}}^j}^\prime(P^{(m)}_+) & \phi_{\mu_{v_{s}}^j}^\prime(P^{(m)}_+)
			                               \end{array}\right].
	\end{aligned}
\]
The matrix \(\bm Q\) formed by SAW basis is given by
\begin{equation}\label{eqn:SAWBasisCoarseMat}
	\bm Q =
	\begin{bmatrix}
		\bm Q^{(1)} &             &        &               \\
		            & \bm Q^{(2)} &        &               \\
		            &             & \ddots &               \\
		            &             &        & \bm Q^{(N-1)}
	\end{bmatrix}.
\end{equation}
We are now ready to introduce the two-level algorithm. Let
\begin{equation}\label{eqn:splitLambda}
	\bm\lambda_k=\tilde{\bm\lambda}_k+\bm Qy_k.
\end{equation}
Then, from \eqref{eqn:Orth2Q} we can deduce that
\[
	\begin{aligned}
		            & \bm Q^*\bm r_k=\bm Q^*(\bm d-\bm F\bm\lambda_k)=0\quad\Rightarrow\quad \bm Q^*\bm F\bm Q\bm y_k=\bm Q^*(\bm d-\bm F\tilde{\bm\lambda}_k) \\
		\Rightarrow & \bm y_k=(\bm Q^*\bm F\bm Q)^{-1}\bm Q^*\bm d-(\bm Q^*\bm F\bm Q)^{-1}\bm Q^*\bm F\tilde{\bm\lambda}_k.
	\end{aligned}
\]
Substituting \(\bm y_k\) into \eqref{eqn:splitLambda}, we obtain
\begin{equation}\label{eqn:defProjector}
	{\bm\lambda}_k=\underbrace{\bm Q(\bm Q^{*}\bm F\bm Q)^{-1}\bm Q^{*}\bm d}_{\bm\lambda_0}+\underbrace{(\bm I-\bm Q(\bm Q^{*}\bm F\bm Q)^{-1}\bm Q^{*}\bm F)}_{\bm P}\tilde{\bm\lambda}_k.
\end{equation}
It is easy to verify that \(\bm P\) is a projection matrix and \(\bm P\bm Q=\bm 0\).

In summary, the solver of \eqref{eqn:LagMultipliersEqn} in the FETI-SAW algorithm is essentially an iterative algorithm that takes \(\bm\lambda_0\) as the initial value to solve
\begin{equation}\label{eqn:precondSystem}
	\bm M^{-1}\bm F\bm\lambda=\bm M^{-1}\bm d,\quad \bm\lambda\in \bm\lambda_0+{\rm Range}(\bm P),
\end{equation}
where \(\bm M^{-1}\) is the Dirichlet preconditioner or the lumped preconditioner.
\begin{remark}
	Typically, a favorable performance can be achieved by constructing the Rayleigh wave coarse space with only 1--3 distinct velocities. Each type of SAW device illustrated in Figure \ref{fig:CommonSAWdevices} can be regarded as a semi-infinite piezoelectric substrate with a partially metallized surface. Thus, incorporating the metallized-surface Rayleigh wave velocity into the selection of SAW velocity typically yields excellent performance. Readers may refer to \cite{Morgan-2007} for detailed procedures to compute the metallized-surface Rayleigh wave velocity.
\end{remark}
\subsection{Key points of implementation}\label{sect:implementation}
For SAW equations, the Dirichlet preconditioner generally outperforms the lumped preconditioner significantly. Thus, the FETI-SAW algorithm adopts the Dirichlet preconditioner by default. Then, we present the implementation details of the FETI-SAW algorithm in the following three stages.

\paragraph{Initialization} We first assemble the matrices corresponding to the five types of subdomains illustrated in Figure \ref{fig:subDomain}, which we denote sequentially as \(\bm K^{L}\), \(\bm K^{r}\), \(\bm K^{b}\), \(\bm K^{t}\), and \(\bm K^{R}\), and generate the right-hand vector for the odd-numbered IDTs. We then initialize the preconditioner by generating the submatrix blocks for each subdomain (i.e., \(\bm K_{i\Gamma}^{\cdot}\), \(\bm K_{\Gamma\Gamma}^{\cdot}\), and the LU decomposition of \(\bm K_{ii}^{\cdot}\)). Next, we perform LU decomposition on the matrices for each type of subdomain.
Finally, we construct the Rayleigh wave basis functions based on the given velocities to form the matrix \(\bm Q\) and perform LU decomposition on \(\bm Q^*\bm F\bm Q\).

\paragraph{Iteration} In this paper, we design the solver for equation \eqref{eqn:LagMultipliersEqn} based on the generalized conjugate residual (GCR) algorithm, and the procedure is presented in Algorithm \ref{alg:FETI-SAW}. In Algorithm \ref{alg:FETI-SAW}, we denote the vector obtained by multiplying matrix \(\bm A\) with vector \(y\) at the \(k\)-th step as \([\bm A y]_k\), the \(l_2\) inner product of vectors \(x\) and \(y\) as \((x,y)\) and \(l^2\) norm as \(\|\cdot\|_2\).

\begin{algorithm}[htbp]
	\caption{two-level GCR algorithm}\label{alg:FETI-SAW}
	\footnotesize
	\begin{algorithmic}[1]
		\State Input: \(epcl\), \(maxit\)
		\State \(\bm\lambda_0=\bm Q(\bm Q^{*}\bm F\bm Q)^{-1}\bm Q^{*}\bm d\)
		\State \(\bm r_0=\bm d-\bm F\bm\lambda_0\), \(\bm z_0=\bm M_D^{-1}\bm r_0\), \(\bm y_0=\bm P\bm z_0\)
		\State \(\bm p_0=\bm y_0\), \([\bm F\bm p]_0=\bm F\bm p_0\), \([\bm M_D^{-1}\bm F\bm p]_{0}=\bm M_D^{-1}[\bm F\bm p]_0\)
		\For{$k = 1, 2, \cdots, maxit$}
		\State \(\alpha_k\)=$\frac{\left(\bm z_{k-1},[\bm M_D^{-1}\bm F\bm p]_{k-1}\right)}{\|[\bm M_D^{-1}\bm F\bm p]_{k-1}\|^2_{2}}$
		\State \(\bm\lambda_{k}=\bm\lambda_{k-1}+\alpha_k\bm p_{k-1}\)
		\State \(\bm r_k=\bm r_{k-1}-\alpha_k[\bm F\bm p]_{k-1}\)
		\If{\(\|\bm r_k\|_{2}/\|\bm d\|_{2}<epcl\)}
		\State Return \(\bm\lambda_{k}\)
		\EndIf
		\State \(\bm z_k=\bm z_{k-1}-\alpha_k [\bm M_D^{-1}\bm F\bm p]_{k-1}\)
		\State \(\bm y_k=\bm P\bm z_k\), \([\bm F\bm y]_k=\bm F\bm y_k\), \([\bm M_D^{-1}\bm F\bm y]_k=\bm M_D^{-1}[\bm F\bm y]_k\)
		\State \(\beta_{kj}=-\frac{\left([\bm M_D^{-1}\bm F\bm y]_k,[\bm M_D^{-1}\bm F\bm p]_{j}\right)}{\|[\bm M_D^{-1}\bm F\bm p]_{j}\|_{2}^2}\), for \(j=0,1,\cdots,k-1\)
		\State \(\bm p_k=\bm y_k+\sum_{j=0}^{k-1}\beta_{kj}\bm p_j\)
		\State \([\bm F\bm p]_k=[\bm F\bm y]_k+\sum_{j=0}^{k-1}\beta_{kj}[\bm F\bm p]_j\)
		\State \([\bm M_D^{-1}\bm F\bm p]_k=[\bm M_D^{-1}\bm F\bm y]_k+\sum_{j=0}^{k-1}\beta_{kj}[\bm M_D^{-1}\bm F\bm p]_j\)
		\EndFor
	\end{algorithmic}
\end{algorithm}
\paragraph{Computation of the original solution} After the iterative process, we obtain the Lagrange multipliers. Then, we can compute the original solution $\bm X$ in parallel by subdomain using the formulation in \eqref{eqn:computeSourceSol}. The algorithm is then executed completely.

As described in the Initialization stage, the use of the Dirichlet preconditioner only requires performing LU decomposition on the matrix $\bm K_{ii}^{\cdot}$ in each type of subdomain. The  cost of this stage is entirely within an acceptable range, as it does not increase with the number of subdomains. When calculating $\bm F\bm Q$ in $\bm Q^*\bm F\bm Q$, since both $\bm F$ and $\bm Q$ adopt a subdomain-based block structure, parallel computation can be implemented. The resulting coarse problem matrix $\bm Q^*\bm F\bm Q$ is a square matrix with the same dimension as the coarse space and has a small size. Given that the matrices within the same type of subdomain are identical, it is unnecessary to assemble the matrix for each subdomain. Instead, we only need to map each subdomain to the corresponding matrix among $\bm K^{L}$, $\bm K^{r}$, $\bm K^{b}$, $\bm K^{t}$, and $\bm K^{R}$ according to the structure of the SAW device.

We briefly analyze the computational cost of Algorithm \ref{alg:FETI-SAW}. In Algorithm \ref{alg:FETI-SAW}, line \(13\) accounts for the major computational cost in each iterative step. In this step, the operation that consumes the most computational resources is solving the linear system, which includes solving one coarse problem, performing two multiplications of matrix \(\bm F\) with vectors, and one multiplication of the Schur complement with vector. We denote the size of \(\bm K^{j}\) (\(j\in\{L,r,b,l,R\}\)) as \(n_{j}\times n_{j}\) and the size of \(\bm Q^*\bm F\bm Q\) as \(n_Q\times n_Q\), where \(n_Q=16N_v(N-1)\). Since the size of \(\bm K^j_{ii}\) varies linearly with \(n_{j}\), we also denote it as \(n_{j}\). Then, in each iteration, the total cost of this step is
\begin{equation}\label{eqn:computionalofIter}
	O\left(3(n_{L}^2+n_R^2+N_rn_r^2+N_bn_b^2+N_tn_t^2)+n_Q^2\right),
\end{equation}
since we have performed the LU decomposition on the relevant matrices. In line \(13\), only the coarse problem is global, while we perform all other operations independently in the subdomains, and \(n_Q\) is typically very small. That is, in \eqref{eqn:computionalofIter}, the operations corresponding to the terms inside the parentheses can be performed in subdomain-parallel,
while only the small-size coarse problem associated with the $n_Q^2$ term needs to be solved serially. Thus, through our analysis, we conclude that this algorithm has good parallelism and low memory resource requirements.
\begin{remark}
	Sometimes, the Rayleigh wave basis constructed from the selected wave velocities exhibits strong linear correlation. Therefore, to ensure the stability of the algorithm, we need to check the basis and eliminate some basis functions, so that the resulting matrix \(\bm Q\) is column-full-rank. In practice, our treatment is as follows: calculate the singular values of each matrix \(\bm Q^{(m)}\), and eliminate the basis functions corresponding to the singular values whose absolute values are smaller than the preset threshold. Obviously, this process can also be parallelized by subdomains, and its size is small, exerting no significant impact on the overall runtime of the algorithm.
\end{remark}
\section{Numerical examples}\label{sect:NumericalExample}
\subsection{Notations and default settings}
To present the numerical examples clearly, we introduce some additional symbols except those illustrated in Figure \ref{fig:2_5DresonatorPML} and some default settings.

For SAW devices, the typical wavelength is generally equal to $2p_t$. Since the SAWs propagate near the surface of the piezoelectric substrate \cite{Morgan-2007}, we fix the length of the piezoelectric substrate along the $x_3$-direction at three typical wavelengths. We denote the electrode lengths along the $x_1$- and $x_3$-directions as $d^e$ and $h^e$, respectively. For ease of presentation, we introduce the frequency \(f\), which relates to the angular frequency \(\omega\) via the relation \(\omega=2\pi f\). For the PML, $d_{x_1}$ is set to one typical wavelength and $d_{x_3}$ to $2/3$ of the typical wavelength. The damping functions $\eta_i(x_i)$ of the PML take the form
\[
	\eta_i(x_i)=\frac{10^6}{\omega}\left(\frac{x_i-\hat{x}_i}{d_{x_i}}\right)^3,\quad i=1,3,
\]
where $\hat{x}_i$ denotes the $x_i$-coordinate of the interface between the PML domain and the piezoelectric substrate.

In the following numerical examples, we use two types of piezoelectric material, 128YXLN and 42YXLT, and denote their metallized-surface Rayleigh wave velocities as $v_{\rm LN}$ and $v_{\rm LT}$, respectively. For both types of the piezoelectric substrates, we construct the Rayleigh wave coarse space based on the set of wave velocities given by
\[
	\bm{v}
	\ =
	\begin{cases}
		\{v_{\rm LN},\; 2000\ \text{m/s}\}, & \text{for } 128\text{YXLN}, \\
		\{v_{\rm LT},\; 2000\ \text{m/s}\}, & \text{for } 42\text{YXLT}.
	\end{cases}
\]

For the discretization process, we adopt the Q2 FEM to discretize each subdomain, with consistent grid points enforced at the interfaces between subdomains and a single element fixed along the \(x_2\)-direction. We denote the number of elements per typical wavelength as \(N_{\rm elem}\).

During the iterative process, we define the relative residual at the \(k\)-th iteration step as
\begin{equation}\label{eqn:relErrLagrangeMulti}
	res_{k}=\frac{\|\bm d-\bm F\bm\lambda_k\|_2}{\|\bm d\|_2},
\end{equation}
and specify the termination residual threshold as \(10^{-6}\). We denote the total number of iterations as \(iter\).  The runtime data reported in the numerical results are the average values of the identical process repeated 10 times.

All numerical examples were implemented in MATLAB R2024a on a high-performance computing cluster node of Tongji University Intelligent Computing Center, equipped with a AMD EPYC 9654 processor, where 64 cores were allocated for computations.
\subsection{IDTs}
For convenient design of controlled experiments, we first adopt the IDT problems illustrated in Figure \ref{fig:CommonSAWdevices} (left) to investigate the effects of variations in the number of subdomains, frequency, and mesh size on the iteration numbers and runtime of the FETI-SAW algorithm. In this subsection, we mainly focus on the performance of the solver for equation \eqref{eqn:LagMultipliersEqn}. We therefore denote the time consumed per iterative step and that for the entire iterative process as \textit{per-step runtime} and \textit{total runtime}, respectively. In this numerical experiment, we default to using $\bm{v}$ to construct the Rayleigh wave coarse space, whose dimension is less than $32(N_t-1)$. The parameters of the model employed in this subsection are listed in Table \ref{tab:idt_params}.

\begin{table}[H]
	\footnotesize
	\centering
	\begin{tabular}{l|ccccc}
		\toprule
		parameter & substrate & electrode & \(p_t\) (\(\mu \rm m\)) & \(d^e\) (\(\mu \rm m\)) & \(h^e\) (\(\mu \rm m\)) \\
		\midrule
		para1     & 128YXLN   & Al        & $1$                     & $0.5$                   & $0.15$                  \\
		para2     & 42YXLT    & Cu        & $1$                     & $0.5$                   & $0.15$                  \\
		\bottomrule
	\end{tabular}
	\caption{Two sets of parameters.}
	\label{tab:idt_params}
\end{table}
Throughout this study, we have explored a variety of algorithm design schemes. To better demonstrate the efficiency and distinctiveness of the FETI-SAW algorithm, we first briefly present the numerical results of several seemingly feasible alternatives.

We denote the standard GCR algorithm for equation \eqref{eqn:LagMultipliersEqn} as Algorithm 2, the variant employing only the Dirichlet preconditioner as Algorithm 3, the variant adopting only the Rayleigh wave coarse space as Algorithm 4, and the variant combining the lumped preconditioner with the Rayleigh wave coarse space as Algorithm 5. Then, we fix $N_{\rm elem}=20$ and solve the IDT problems with parameter $\text{para1}$ using the above four algorithms for different values of $N_t$ and frequencies. The numerical results are shown in Table \ref{tab:failed idea}
\begin{table}[H]
	\scriptsize
	\centering
	\begin{tabular}{l|lllllll}
		\toprule
		\multirow{3}{*}{Algorithm}
		 & \(f \)
		 & \multicolumn{3}{c}{100 MHz}
		 & \multicolumn{3}{c}{200 MHz}                   \\
		\cmidrule(lr){2-2} \cmidrule(lr){3-5} \cmidrule(lr){6-8}
		 & \(N_t\)
		 & 11                          & 21     & 31
		 & 11                          & 21     & 31     \\
		\midrule
		\multirow{2}{*}{2}
		 & iter
		 & 523                         & 726    & 895
		 & 502                         & 705    & 874    \\
		 & total runtime (s)
		 & 88.01                       & 175.34 & 246.31
		 & 84.02                       & 169.50 & 240.85 \\
		\addlinespace[1ex]
		\multirow{2}{*}{3}
		 & iter
		 & 62                          & 94     & 126
		 & 84                          & 140    & 197    \\
		 & total runtime (s)
		 & 19.41                       & 39.74  & 59.07
		 & 26.44                       & 59.92  & 92.70  \\
		\addlinespace[1ex]
		\multirow{2}{*}{4}
		 & iter
		 & 356                         & 442    & 501
		 & 308                         & 372    & 413    \\
		 & total runtime (s)
		 & 115.49                      & 196.85 & 367.53
		 & 101.33                      & 165.31 & 214.21 \\
		\addlinespace[1ex]
		\multirow{2}{*}{5}
		 & iter
		 & 55                          & 63     & 69
		 & 46                          & 51     & 53     \\
		 & total runtime (s)
		 & 17.43                       & 28.34  & 37.00
		 & 14.75                       & 23.29  & 27.80  \\
		\bottomrule
	\end{tabular}
	\caption{Numerical results of the four underperforming algorithms.}
	\label{tab:failed idea}
\end{table}

It can be seen from Table \ref{tab:failed idea} that the standard GCR algorithm (Algorithm 2) requires a considerable number of iterations to converge. After only introducing the coarse space to construct the two-level algorithm (Algorithm 4), the number of iterations of the algorithm decreases significantly but remains relatively large. When only the Dirichlet preconditioner is used (Algorithm 3), the number of iterations is greatly reduced, yet the algorithm is still sensitive to the changes in the number of subdomains and frequencies. After combining the lumped preconditioner with the Rayleigh wave coarse space (Algorithm 5), we observe a substantial reduction in the algorithm's iteration count, as well as reduced sensitivity to subdomain number and frequency. However, such efficiency still cannot qualify the algorithm as a fast solver for large-size simulation.

Next, we solve the IDT problems with $\text{para1}$ and $\text{para2}$ as parameters, respectively, using the FETI-SAW algorithm for different values of $N_t$, $f$ and $N_{\rm elem}$. The numerical results are illustrated in Tables \ref{tab:FETISAWIDTIter}--\ref{tab:FETISAWIDTItervarinElem}.

\begin{table}[H]
	\scriptsize
	\centering
	\setlength{\tabcolsep}{5pt}
	\begin{tabular}{c|cccc|cccc}
		\toprule
		\multicolumn{1}{c|}{parameter}
		 & \multicolumn{4}{c|}{para1}
		 & \multicolumn{4}{c}{para2}                   \\
		\midrule
		\diagbox{$f$}{\(iter\)}{$N_t$}
		 & 21                         & 81 & 141 & 201
		 & 21                         & 81 & 141 & 201 \\
		\midrule
		100 MHz
		 & 9                          & 9  & 9   & 9
		 & 12                         & 12 & 12  & 12  \\
		200 MHz
		 & 7                          & 7  & 7   & 8
		 & 4                          & 3  & 3   & 3   \\
		300 MHz
		 & 6                          & 6  & 7   & 8
		 & 7                          & 10 & 11  & 8   \\
		400 MHz
		 & 6                          & 8  & 8   & 8
		 & 8                          & 9  & 14  & 11  \\
		\bottomrule
	\end{tabular}

	\caption{Iteration numbers of the FETI-SAW algorithm (\(N_{\rm elem}=20\)).}
	\label{tab:FETISAWIDTIter}
\end{table}
\begin{table}[H]
	\fontsize{7pt}{7.2pt}\selectfont
	\centering
	\setlength{\tabcolsep}{2pt}
	\begin{tabular}{c|cccc|cccc}
		\toprule
		\multicolumn{1}{c|}{parameter}
		 & \multicolumn{4}{c|}{para1}
		 & \multicolumn{4}{c}{para2}                                              \\
		\midrule
		\diagbox{$f$}{runtime (s)}{$N_{t}$}
		 & 21                         & 81          & 141          & 201
		 & 21                         & 81          & 141          & 201          \\
		\midrule
		100 MHz
		 & 4.73 (0.53)                & 7.78 (0.86) & 13.87 (1.54) & 22.34 (2.48)
		 & 5.95 (0.50)                & 9.69 (0.81) & 16.96 (1.41) & 27.22 (2.27) \\
		200 MHz
		 & 3.74 (0.53)                & 6.43 (0.92) & 12.39 (1.77) & 22.51 (2.81)
		 & 2.23 (0.56)                & 3.24 (1.08) & 6.48 (2.16)  & 10.75 (3.58) \\
		300 MHz
		 & 3.34 (0.56)                & 6.24 (1.04) & 14.54 (2.08) & 27.44 (3.43)
		 & 3.64 (0.52)                & 9.33 (0.93) & 20.19 (1.84) & 24.93 (3.12) \\
		400 MHz
		 & 3.28 (0.55)                & 8.08 (1.01) & 16.43 (2.05) & 27.65 (3.46)
		 & 4.13 (0.52)                & 8.44 (0.94) & 25.18 (1.80) & 33.36 (3.03) \\
		\bottomrule
	\end{tabular}
	\caption{Runtime of the FETI-SAW algorithm (\(N_{\rm elem}=20\); data outside the bracket is total runtime;  data inside is per-step runtime).}
	\label{tab:FETISAWIDTRuntime}
\end{table}
\begin{table}[H]
	\scriptsize
	\centering
	\begin{tabular}{c|ccccc|ccccc}
		\toprule
		\multicolumn{1}{c|}{parameter}
		 & \multicolumn{5}{c|}{para1}
		 & \multicolumn{5}{c}{para2}                      \\
		\midrule
		\diagbox{$f$}{\(iter\)}{$N_{\rm elem}$}
		 & 20                         & 25 & 30 & 35 & 40
		 & 20                         & 25 & 30 & 35 & 40 \\
		\midrule
		100 MHz
		 & 9                          & 9  & 9  & 9  & 9
		 & 12                         & 12 & 12 & 12 & 12 \\
		200 MHz
		 & 7                          & 7  & 7  & 7  & 7
		 & 4                          & 3  & 3  & 3  & 3  \\
		300 MHz
		 & 6                          & 6  & 6  & 6  & 6
		 & 7                          & 7  & 7  & 7  & 7  \\
		400 MHz
		 & 6                          & 6  & 6  & 6  & 6
		 & 8                          & 8  & 8  & 8  & 8  \\
		\bottomrule
	\end{tabular}
	\caption{Iteration numbers of the FETI-SAW algorithm (\(N_t=21\)).}
	\label{tab:FETISAWIDTItervarinElem}
\end{table}
In Tables \ref{tab:FETISAWIDTIter} and \ref{tab:FETISAWIDTRuntime}, we fixed \(20\) elements per typical wavelength while varying the number of subdomains (\(N_t\)) and the frequency (\(f\)). Table \ref{tab:FETISAWIDTIter} shows that the iteration numbers remain almost unchanged as $N_t$ increases for a fixed frequency. For a fixed number of subdomains, the iteration counts of both models remain stable as the frequency varies. Table \ref{tab:FETISAWIDTRuntime} shows that the time consumed per iterative step increases slowly as the number of subdomains increases. When $N_t > 81$, the growth rate of the per-step iteration time is significantly higher than that when $N_t < 81$, which is attributed to the limitation of the CPU core count. Table \ref{tab:FETISAWIDTItervarinElem} shows that the iteration count is independent of \(N_{\text{elem}}\). Compared with the results in Table \ref{tab:failed idea}, the FETI-SAW algorithm exhibits an absolute advantage in solving the IDT problems. Remarkably, we achieve such favorable results simply by combining a low-dimensional coarse space with the Dirichlet preconditioner. Therefore, it does not introduce excessive computational overhead in each iterative step.

In summary, the numerical results in this subsection demonstrate that the FETI-SAW algorithm exhibits insensitivity to the number of subdomains, mesh size, and frequency for solving IDT problems. Furthermore, they also highlight the potential of the FETI-SAW algorithm for addressing large-size SAW device problems.
\subsection{SAW resonators}
In this subsection, we investigate the numerical performance of the FETI-SAW algorithm for solving the SAW resonator problems introduced in section \ref{sect:MathematicalModel}. The two sets of parameters we considered are given in Table \ref{tab:resonator_params}. In practice, we are more concerned with the overall efficiency and accuracy of the algorithm in solving the equation \eqref{eqn:originalEqn}. Therefore, we present the runtime for each stage mentioned in section \ref{sect:implementation} as well as the total runtime, and define the relative residual of the original equation as
\[
	res_X=\frac{\|\bm f-\bm K\bm X_{iter}\|_2}{\|\bm f\|_2}.
\]
\begin{table}[H]
	\footnotesize
	\centering
	\begin{tabular}{l|lllllll}
		\toprule
		parameter & substrate & electrode & \(p_r\) (\(\mu\rm m\)) & \(p_b\) (\(\mu\rm m\)) & \(p_t\) (\(\mu\rm m\)) & \(d^e\) (\(\mu\rm m\)) & \(h^e\) (\(\mu\rm m\)) \\ \midrule
		para1     & 128YXLN   & Cu        & 1                      & 0.5                    & 1                      & 0.5                    & 0.17                   \\
		para2     & 42YXLT    & Al        & 1                      & 0.5                    & 1                      & 0.55                   & 0.16                   \\ \bottomrule
	\end{tabular}
	\caption{Parameters of the SAW resonators.}
	\label{tab:resonator_params}
\end{table}

\begin{figure}[H]
	\centering
	\begin{subfigure}{\linewidth}
		\centering
		\includegraphics[width=\linewidth]{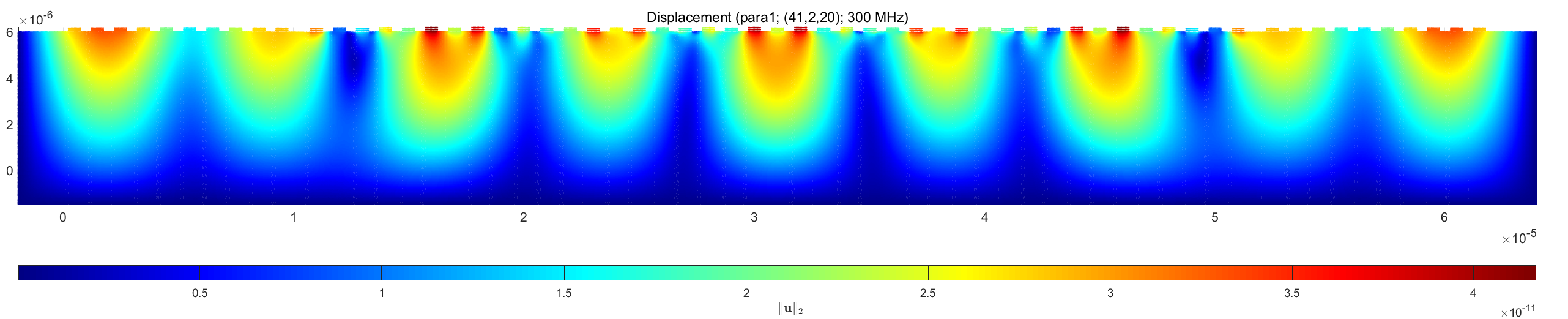}
		\caption{} 
	\end{subfigure}
	\vspace{1em} 

	\begin{subfigure}{\linewidth}
		\centering
		\includegraphics[width=\linewidth]{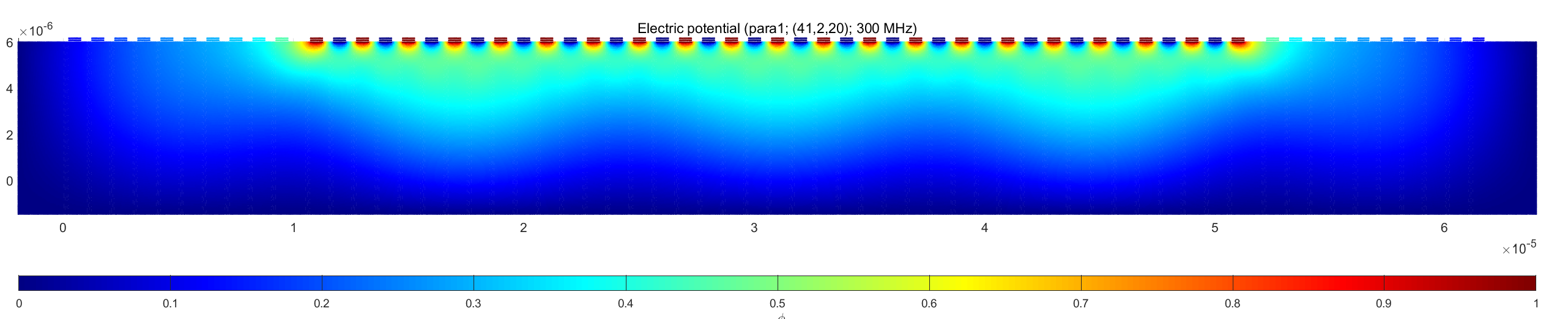}
		\caption{} 
	\end{subfigure}

	\caption{Numerical results of the model with para1, structure $(41,2,20)$ at 300 MHz: (a) displacement; (b) electric potential.}
	\label{fig:para1Struct1}
\end{figure}

\begin{figure}[H]
	\centering
	\begin{subfigure}{\linewidth}
		\centering
		\includegraphics[width=\linewidth]{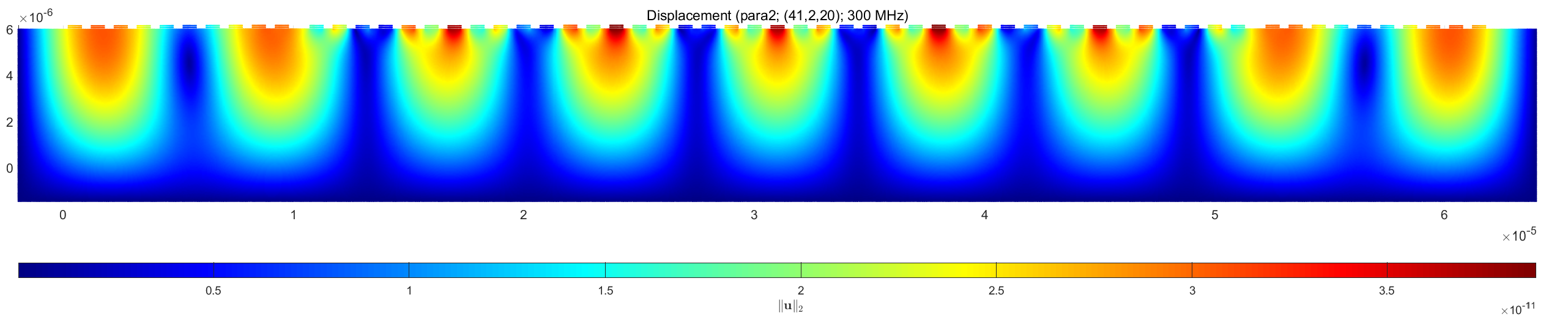}
		\caption{} 
	\end{subfigure}
	\vspace{1em} 

	\begin{subfigure}{\linewidth}
		\centering
		\includegraphics[width=\linewidth]{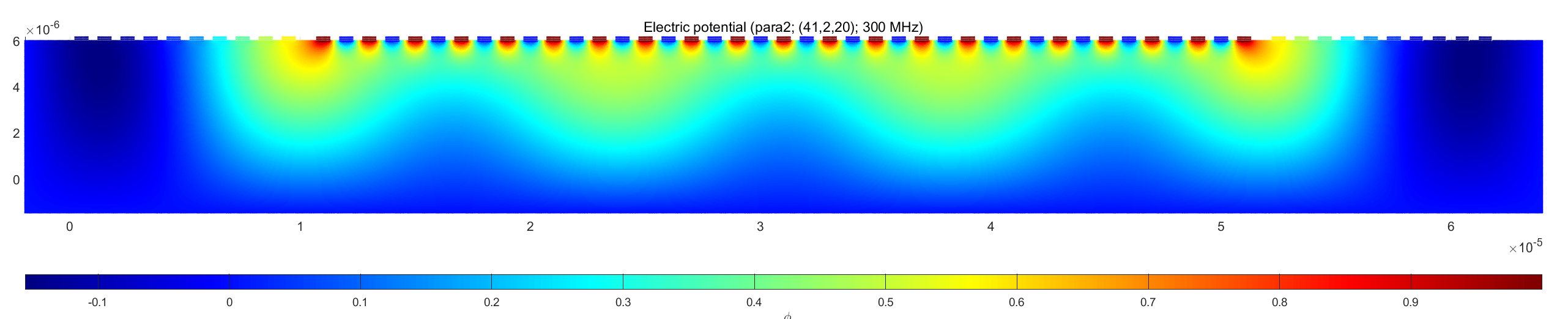}
		\caption{} 
	\end{subfigure}

	\caption{Numerical results of the model with para2, structure $(41,2,20)$ at 300 MHz: (a) displacement; (b) electric potential.}
	\label{fig:para2Struct1}
\end{figure}

\begin{figure}[H]
	\centering
	\begin{subfigure}{\linewidth}
		\centering
		\includegraphics[width=\linewidth]{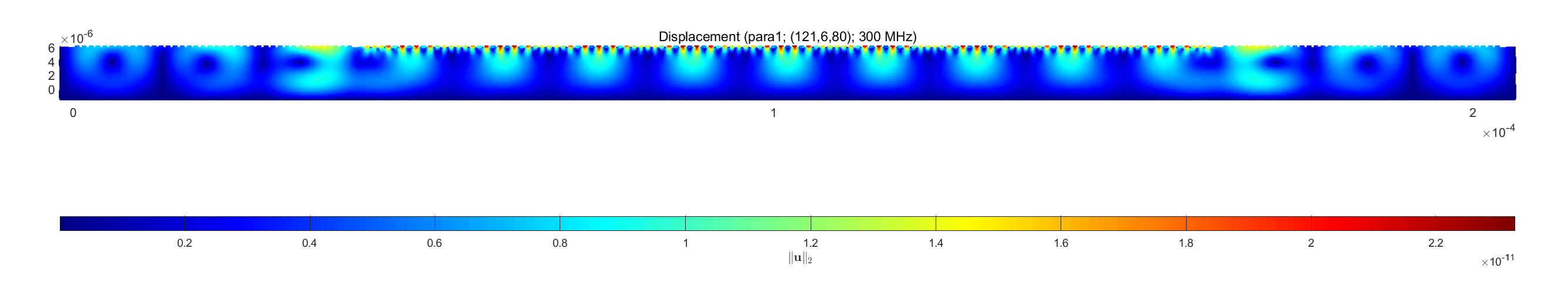}
		\caption{} 
	\end{subfigure}
	\vspace{1em} 

	\begin{subfigure}{\linewidth}
		\centering
		\includegraphics[width=\linewidth]{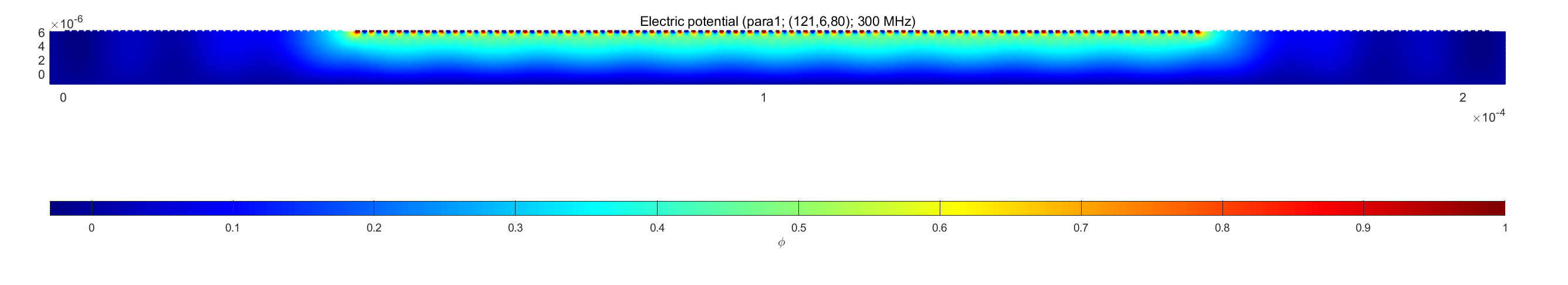}
		\caption{} 
	\end{subfigure}

	\caption{Numerical results of the model with para1, structure $(121,6,80)$ at 300 MHz: (a) displacement; (b) electric potential.}
	\label{fig:para1}
\end{figure}

\begin{figure}[H]
	\centering
	\begin{subfigure}{\linewidth}
		\centering
		\includegraphics[width=\linewidth]{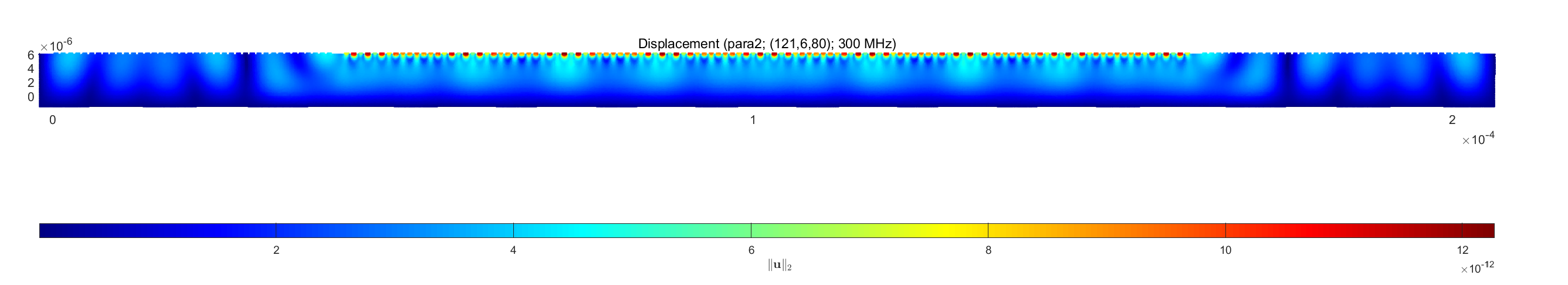}
		\caption{}
	\end{subfigure}
	\vspace{1em}

	\begin{subfigure}{\linewidth}
		\centering
		\includegraphics[width=\linewidth]{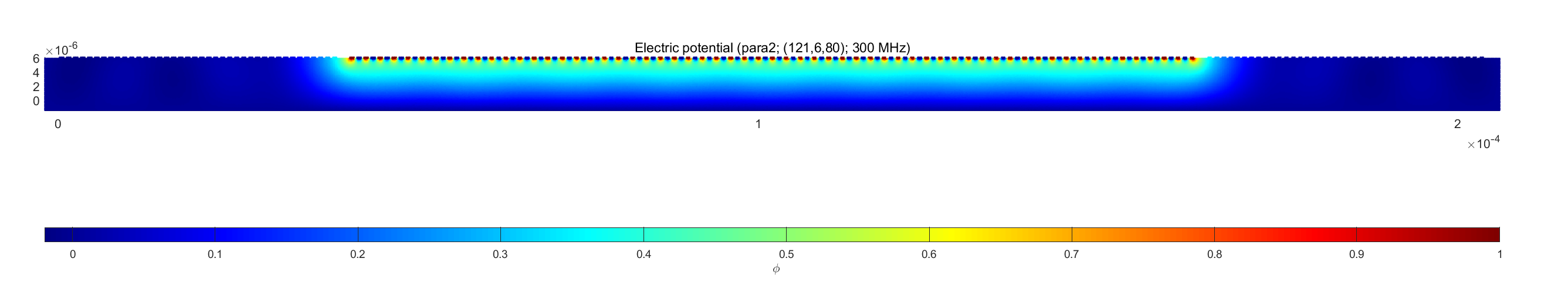}
		\caption{}
	\end{subfigure}

	\caption{Numerical results of the model with para2, structure $(121,6,80)$ at 300 MHz: (a) displacement; (b) electric potential.}
	\label{fig:para2}
\end{figure}

We fix $N_{\rm elem}=20$ and apply the FETI-SAW algorithm to solve several sets of resonator problems with the two parameter sets presented in Table \ref{tab:resonator_params}. The corresponding numerical results are summarized in Tables \ref{tab:FETISAWResonatorPara1}--\ref{tab:FETISAWResonatorPara2} , and the displacement and electric potential distributions for part of these models are illustrated in Figure. \ref{fig:para1Struct1} -- \ref{fig:para2}.

\begin{table}[H]
	\scriptsize
	\centering
	\begin{tabular}{llllllllll}
		\toprule
		\multicolumn{10}{c}{para1}                                                                                                                                                                                                  \\ \midrule
		\multicolumn{2}{l|}{(Nt,Nb,Nr)}                  & \multicolumn{4}{c|}{(41,2,20)}      & \multicolumn{4}{c}{(121,6,80)}                                                                                                     \\ \midrule
		\multicolumn{2}{l|}{$f$ (MHz)}                   & 100                                 & 200                            & 300   & \multicolumn{1}{l|}{400}  & 100                        & 200   & 300    & 400             \\ \midrule
		\multicolumn{2}{l|}{iter}                        & 14                                  & 11                             & 10    & \multicolumn{1}{l|}{10}   & 19                         & 14    & 12     & 18              \\
		\multicolumn{2}{l|}{$res_X (10^{-13})$}          & 0.16                                & 3.38                           & 0.34  & \multicolumn{1}{l|}{0.46} & 0.15                       & 3.53  & 0.31   & 8.03            \\
		\multicolumn{1}{c}{\multirow{4}{*}{runtime (s)}} & \multicolumn{1}{l|}{initialization} & 9.28                           & 10.53 & 14.38                     & \multicolumn{1}{l|}{14.32} & 45.75 & 60.58  & 99.68  & 99.52  \\
		\multicolumn{1}{c}{}                             & \multicolumn{1}{l|}{iteration}      & 11.43                          & 9.16  & 9.00                      & \multicolumn{1}{l|}{8.97}  & 47.30 & 40.39  & 41.44  & 63.59  \\
		\multicolumn{1}{c}{}                             & \multicolumn{1}{l|}{originalSol}    & 0.13                           & 0.13  & 0.13                      & \multicolumn{1}{l|}{0.13}  & 0.41  & 0.40   & 0.41   & 0.40   \\
		\multicolumn{1}{c}{}                             & \multicolumn{1}{l|}{total}          & 20.83                          & 19.82 & 23.51                     & \multicolumn{1}{l|}{23.43} & 93.46 & 101.37 & 141.53 & 163.51 \\
		\multicolumn{2}{l|}{DOFs $(10^6)$}               & 1.63                                & 1.63                           & 1.63  & \multicolumn{1}{l|}{1.63} & 5.15                       & 5.15  & 5.15   & 5.15            \\
		\multicolumn{2}{l|}{dim($\bm Q$)}                & 1024                                & 1280                           & 1792  & \multicolumn{1}{l|}{1792} & 3328                       & 4160  & 5824   & 5824            \\ \bottomrule
	\end{tabular}
	\caption{Numerical results of the SAW resonator with parameter para1.}
	\label{tab:FETISAWResonatorPara1}
\end{table}

\begin{table}[H]
	\scriptsize
	\centering
	\begin{tabular}{llllllllll}
		\toprule
		\multicolumn{10}{c}{para2}                                                                                                                                                                                                  \\ \midrule
		\multicolumn{2}{l|}{(Nt,Nb,Nr)}                  & \multicolumn{4}{c|}{(41,2,20)}      & \multicolumn{4}{c}{(121,6,80)}                                                                                                     \\ \midrule
		\multicolumn{2}{l|}{$f$ (MHz)}                   & 100                                 & 200                            & 300   & \multicolumn{1}{l|}{400}  & 100                        & 200   & 300    & 400             \\ \midrule
		\multicolumn{2}{l|}{iter}                        & 15                                  & 8                              & 10    & \multicolumn{1}{l|}{13}   & 22                         & 8     & 13     & 26              \\
		\multicolumn{2}{l|}{$res_X (10^{-13})$}          & 0.20                                & 0.55                           & 1.63  & \multicolumn{1}{l|}{0.13} & 0.21                       & 0.55  & 0.21   & 0.11            \\
		\multicolumn{1}{c}{\multirow{4}{*}{runtime (s)}} & \multicolumn{1}{l|}{initialization} & 8.59                           & 12.43 & 12.46                     & \multicolumn{1}{l|}{12.42} & 37.09 & 78.70  & 78.15  & 78.92  \\
		\multicolumn{1}{c}{}                             & \multicolumn{1}{l|}{iteration}      & 11.64                          & 6.99  & 8.58                      & \multicolumn{1}{l|}{11.09} & 51.50 & 26.70  & 40.86  & 80.68  \\
		\multicolumn{1}{c}{}                             & \multicolumn{1}{l|}{originalSol}    & 0.14                           & 0.13  & 0.13                      & \multicolumn{1}{l|}{0.13}  & 0.41  & 0.40   & 0.40   & 0.40   \\
		\multicolumn{1}{c}{}                             & \multicolumn{1}{l|}{total}          & 20.36                          & 19.55 & 21.17                     & \multicolumn{1}{l|}{23.63} & 89.00 & 105.80 & 119.41 & 160.01 \\
		\multicolumn{2}{l|}{DOFs $(10^6)$}               & 1.63                                & 1.63                           & 1.63  & \multicolumn{1}{l|}{1.63} & 5.15                       & 5.15  & 5.15   & 5.15            \\
		\multicolumn{2}{l|}{dim($\bm Q$)}                & 912                                 & 1536                           & 1536  & \multicolumn{1}{l|}{1536} & 2940                       & 4992  & 4992   & 4992            \\ \bottomrule
	\end{tabular}
	\caption{Numerical results of the SAW resonator with parameter para2.}
	\label{tab:FETISAWResonatorPara2}
\end{table}
In Tables \ref{tab:FETISAWResonatorPara1} and \ref{tab:FETISAWResonatorPara2}, DOFs denotes the number of DOFs in the original equation, and $\dim(\bm Q)$ represents the dimension of the Rayleigh wave coarse space. Since the structure of the SAW resonators and the waves propagating within them are more complex than those of the IDT problems, the performance of the FETI-SAW algorithm degrades slightly, yet it remains highly efficient. The FETI-SAW algorithm solves the resonator problem with 5 million DOFs in less than three minutes. Furthermore, the relative residual of the original problem reaches \(10^{-13}\), even though we terminate the iterative process when \(res_{\text{iter}}<10^{-6}\), which indicates that the solutions inside each subdomain have achieved sufficiently high accuracy when the iteration terminates.  In summary, the FETI-SAW algorithm exhibits favorable efficiency and accuracy in SAW resonator problems.

\section{Conclusion}
In this paper, we propose an efficient domain decomposition algorithm, referred to as the FETI-SAW algorithm, for the FE simulation of SAW devices. The proposed method exhibits good scalability and efficiency. In the FETI-SAW algorithm, the key to achieving such satisfactory performance lies in the proposed Rayleigh wave coarse space, which is based on the characteristics of the wave field in SAW devices. This coarse space is further combined with the Dirichlet preconditioner to improve the iterative algorithm. In Section \ref{sect:implementation}, we detail the key procedures of the FETI-SAW algorithm and demonstrate its favorable parallel performance. Numerical experiments in Section \ref{sect:NumericalExample} demonstrate the efficiency and scalability of the FETI-SAW algorithm for the FE simulation of IDTs and SAW resonators.

As typical examples, we demonstrate the performance of the FETI-SAW algorithm for the devices including IDTs and SAW resonators. In fact, given the similar physical principles of SAW devices, the proposed algorithm can be directly applied to the simulation of some other SAW devices.
For the reader's convenience, the MATLAB source code accompanying this paper is available at \url{https://github.com/sduwhqlk1997/FETI}, and the C++ version is available at \url{https://github.com/sduwhqlk1997/FETI-SAW-CPP}. These codes can be readily applied to SAW problems with other structures.
\section{Acknowledgements}
The computations in this paper were supported by Tongji University Intelligent Computing Center.

%





\bibliographystyle{plain}
\bibliography{reference}
\end{document}